\documentclass[12pt]{article}
\usepackage{amsfonts}
\usepackage{amssymb}
\usepackage{mathrsfs}
\usepackage{amsmath}

\def\cL{\mathcal L}
\def\cT{\mathcal T}
\def\cD{\mathcal D}
\def\cJ{\mathcal J}
\def\cG{\mathcal G}

\def\cU{\mathcal U}
\def\C{\mathbb{C}}
\def\Z{\mathbb{Z}}
\def\N{\mathbb{N}}

\def\tA{\tilde A}

\newtheorem{theorem}{Theorem}[section]
\newtheorem{prop}{Proposition}[section]

\begin{document}
\centerline{\textbf{\large{A class of $q$-orthogonal polynomial sequences}}}

 \centerline{\textbf{\large{ that extends the $q$-Askey scheme}}}

\bigskip
\centerline{\textbf{Luis Verde-Star}}

\medskip

\centerline{Department of Mathematics, Universidad Aut\'onoma Metropolitana,}

\centerline{ Iztapalapa, Apartado 55-534, M\'exico D. F. 09340, M\'exico.}

\centerline{ E-mail: verde@xanum.uam.mx}

\bigskip 

\begin{abstract}   
	We obtain new explicit formulas for the recurrence coefficients of the $q$-ortho\-go\-nal polynomial sequences in the extended $q$-Hahn class introduced in our previous paper  \cite{Rec}. Our  new formulas  express the recurrence coefficients as rational functions of $q^k$ with numerators and denominators that are completely factored and whose zeroes and poles depend on only four parameters. 
By direct substitution of  particular values for the four parameters we obtain all the polynomial sequences in the $q$-Askey scheme, including the Askey-Wilson and the Racah polynomials. We also obtain some families of sequences that are not mentioned in the book \cite{Hyp}, which is the standard reference for the $q$-Askey scheme.


\end{abstract}

\section{Introduction}
Sequences of orthogonal polynomials that satisfy $q$-difference equations have been studied for a long time. Several classes of such sequences have been studied using diverse approaches. The classes of orthogonal sequences such that their derivatives, differences, or $q$-differences are also orthogonal are called Hahn classes, see \cite{H}. There are numerous recent papers that deal with diverse aspects of classes of $q$-orthogonal polynomials, for example, \cite{Alv}, \cite{KS}, \cite{Koor}, \cite{MBP}, and \cite{Med}.

In the present  paper,  we extend some of  the results obtained in our previous papers \cite{Mops} and \cite{Rec} about  characterizations and explicit formulas for the recurrence coefficients of all the classical discrete orthogonal and $q$-orthogonal polynomial sequences. We present   
 new explicit formulas for the recurrence coefficients of all the  $q$-orthogonal polynomial sequences in a class that extends the $q$-Hahn class and contains all the sequences in the $q$-Askey scheme which are listed  in the book \cite{Hyp}.  We also find some families of $q$-orthogonal sequences that are not mentioned in \cite{Hyp} and may be new.
 
Our new formulas express the recurrence coefficients as rational functions of $q^k$ whose numerators and denominators are written as  products of degree one polynomials in $q^k$. 
If  $\{p_k\}_{k\geqslant 0}$ is  a sequence of monic orthogonal  polynomials we write the corresponding  3-term recurrence relation   in the form
 $$\alpha_k p_{k-1}(x)+\beta_k p_k(x)+p_{k+1}(x)= xp_k(x), \qquad k \geqslant 1, $$
 and define  the sums $\sigma_k= \beta_0 +\beta_1 +\ldots +\beta_k,$ for $k \geqslant 0.$ 
 We show that $\{p_k\}_{k\geqslant 0}$ is in  the extended $q$-Hahn class  if  
  the sequence $\alpha_k$ satisfies  equation (4.3) and  $\sigma_k$ satisfies either (4.4) of (4.8).  Such equations give the coefficients in terms of four parameters that determine the zeroes and poles of $\alpha_k$ as a rational function of $q^k$.  
By direct substitution of suitable values of the parameters we obtain  the recurrence coefficients of all the sequences in the $q$-Askey scheme and many others that may be new. 
For the discrete orthogonal polynomials   
we have obtained analogous results that will appear elsewhere.

  Regarding the  $q$-difference equations satisfied by the elements of the extended $q$-Hahn class,  we found in \cite{Rec} that, modulo a tridiagonal change of basis in the space of polynomials, they are the usual second order equations satisfied by the elements of the classical $q$-Hahn class.

\section{Preliminary material}

  We present next  some definitions and basic properties of  lower semi-matrices  and matrices of orthogonal polynomial sequences that  will be used in the rest of the paper.   For a more detailed account of the theory  see   \cite{Mops}.

We consider  matrices $A=[a_{j,k}]$ over the complex numbers where the indices run over the nonnegative integers. 
 We say that $A$ is a {\sl lower semi-matrix}
 if  there is an integer $m$ such that  $a_{j,k}=0$ whenever
 $j-k < m$.  We denote by $\cL$ the set of all lower semi-matrices.
The entry $a_{j,k}$ lies in the diagonal of index  $n$, (also called  the $n$-th diagonal), if $j-k=n$.
 Therefore, if $m > n$ then the $m$-th diagonal lies below (to the left of) the $n$-th diagonal.
  A nonzero element of $\cL$ is  a diagonal matrix if all of its nonzero elements lie in a single diagonal. 
 Given a nonzero  $A$ in  $\cL$  let $m$ be the minimum integer such that $A$ has at least one nonzero entry in the $m$-th diagonal, then we say that $A$ has {\sl index} $m$ and write ind$(A)=m$. We define the index of the zero matrix to be infinity.

It is clear that the space  $\cL$ is  a complex  vector space with the natural addition of matrices and multiplication by scalars. It is also closed under matrix multiplication, and  if  $A$ and $B$ are  in $\cL$, with ind$(A)=m$ and ind$(B)=n$. Then the product $C=AB$ is a well defined element of $\cL$  and ind$(AB)\geqslant m +n$.
A nonzero $A \in \cL$ of index $m$ is called monic if all of its entries in the $m$-th diagonal are equal to one.

If $A \in \cL$ then  a sufficient (but not necessary) condition for $A$ to have a two-sided inverse is that ind$(A)=0$ and  $a_{k,k}\ne 0$ for $k \geqslant 0$. We denote by $\cG$ the set of all matrices that satisfy such condition. It is easy to see that $\cG$ is a group under matrix multiplication. The unit is the identity matrix $I$ whose entries on the 0-th diagonal are equal to 1 and all other entries are zero.

A matrix $A$  is {\sl banded} if there exists  a pair of integers $(n,m)$ such that $n \leqslant m$ and all the nonzero entries of $A$ lie between the diagonals of indices   $n$ and $m$. In such  case we say that $A$ is $(n,m)$-banded.
The set of banded matrices is closed under addition and multiplication, but a banded matrix may have an inverse that is not  banded.
We denote by $\cT$ the set of $(-1,1)$-banded matrices whose entries in the diagonals of indices  1 and -1 are all nonzero.

A sequence  $\{p_k\}_{k\geqslant 0}$  of polynomials in $\C[x]$  such that
$\deg (p_k(x) - x^k)<k$,  for $k \geqslant 0$,  is said to be a monic polynomial sequence. A  basis  $\cU =\{u_k(x), k \geqslant 0 \}$ of $\C[x]$ such that $\{u_k\}_{k\geqslant 0}$  is a monic polynomial sequence is called triangular basis. We  express  the polynomials $p_k$  of a monic sequence in the form
$$p_k(x)=\sum_{j=0}^k a_{k,j} u_j(x), \qquad a_{k,k}=1, \quad k \geqslant 0,$$
where $\cU$ is a triangular basis. Then we use  the coefficients $a_{k,j}$ to form the  matrix $A=[a_{k,j}]$, 
called  the matrix associated with the sequence    $\{p_k\}_{k\geqslant 0}$  with respect to the basis $\cU$.
Note that   $A$ is a monic element of the group $\cG$. 

A polynomial sequence $\{p_k\}_{k\geqslant 0}$ is {\sl orthogonal} with respect to a linear functional $\mu$ if and only if
$$\mu (p_k p_n)= \gamma_k \delta_{k,n}, \qquad k,n \in \N, \eqno(2.1)$$
 where the $\gamma_k$ are nonzero real or complex  constants.
It is well-known that an orthogonal polynomial sequence satisfies a three-term recurrence relation of the form  
$$\alpha_k p_{k-1}(x)+\beta_k p_k(x)+p_{k+1}(x)= xp_k(x), \qquad k \geqslant 1. \eqno(2.2)$$ 
Such recurrence relation is equivalent to the matrix equation $ L A = A X$, where 
$$L=\left[  \begin{matrix} \beta_0 & 1 & 0 & 0 & \ldots \cr
                     \alpha_1 & \beta_1 & 1 & 0 & \ldots \cr
                     0 & \alpha_2 & \beta_2 & 1 & \ldots \cr
       \vdots &\vdots & \vdots & \vdots & \ddots \cr \end{matrix}
    \right],  \eqno(2.3)$$
 and $X$ is the matrix representation with respect to the basis $\cU$ of the multiplication map  on $\C[x]$ that sends $f(x)$  to $x f(x)$. Note that $L$ is in $\cT$.

 If $A$ is the associated matrix of a monic orthogonal polynomial sequence $\{p_k\}_{k\geqslant 0}$  with respect to some triangular basis, we say that $A$ is a  (monic) {\it matrix of orthogonal polynomials} (abbreviated to {\it MOP}).

\section{The extended Hahn's class of $q$-orthogonal polynomial sequences}
In our previous paper \cite{Rec}  we defined the extended Hahn class of $q$-orthogonal polynomial sequences and we found  explicit formulas for the recurrence coefficients of all of its elements using the initial coefficients $\alpha_1, \beta_0, \beta_1, \beta_2$ and $t$. We include here the definitions and results from  \cite{Rec} that we need in this paper.

 Let $q$ be a nonzero complex number that is not a root of 1. Define  the $q$-difference operator  
$$ \cD_q p(x)= \frac{p(qx) -p(x)}{qx -x}, \qquad p \in \C[x], \eqno(3.1)$$
and the $q$ numbers 
$$  [k] = \frac{q^k -1}{q-1}, \qquad k \in \Z.\eqno(3.2)$$
 Note that  $[k]= 1 + q + q^2+\cdots + q^{k-1}$ if  $k \geqslant 1$ and that  
 $\cD_q x^k= [k] x^{k-1}$   for $ k \geqslant 0. $

 The matrix representation of $\cD_q$  with respect to the basis of  monomials $\{x^n \}_{n\geqslant 0}$, which we denote  by $D_q$, is given by $(D_q)_{k+1,k}= [k+1]$ for $k \geqslant 0$, and all other entries are zero.
 Let $X$ be the matrix representation of the multiplication operator that sends  $p(x)$ to $ x \, p(x)$ with respect to the same  basis  of monomials. Then we have
$$ D_q =\left[ \begin{matrix}  0\  & 0\  & 0\  & 0\  & \ldots \cr [1]\  & 0\  & 0\  & 0\  & \ldots \cr 0\  & [2]\  & 0 \  & 0\  & \ldots \cr 0\  & 0 & [3]\  & 0\  & \ldots \cr \vdots \  & \vdots \  & \vdots \ & \ddots \  & \ \ddots \ \end{matrix}\right],  \qquad  X =   \left[ \begin{matrix}  0\  & 1\  & 0\  & 0\  & 0\ &  \ldots \cr  0\  & 0\  & 1\  & 0\  & 0\ & \ldots \cr  0\  & 0\ &  0\  & 1\  & 0\  & \ldots \cr 0\  & 0\ & 0\ &   0\  & 1\    & \ldots \cr \vdots \  & \vdots \ & \vdots \  & \vdots \ & \ddots \  & \ \ddots \ \end{matrix}\right].    \eqno(3.3)$$

The matrices $D_q$ and $X$ satisfy the identity
$$    X D_q  - q D_q X = I.  \eqno(3.4) $$
The matrix
$$ \hat D_q = \left[ \begin{matrix}  0\  & 1\  & 0\  & 0\  & \ldots \cr 0\  & 0\  & 1/[2]\  & 0\  & \ldots \cr 0\  & 0\  & 0 \  & 1/[3]\  & \ldots \cr 0\  & 0 & 0\  & 0\  & \ddots \cr \vdots \  & \vdots \  & \vdots \ & \vdots \  & \ \ddots \ \end{matrix}\right],  \eqno(3.5)$$
	is a left-inverse of $D_q$, that is $\hat D_q  D_q = I$. 
  
 The Hahn's class $\cJ_q$ of $\cD_q$ is by definition  the set of all MOPs $A$ such that $\tA=\hat D_q A D_q$ is also a MOP. 

Let $A$ be an element of $\cJ_q$ and let $L$ and $M$ be the elements of $\cT$ such that $LA=AX$ and $M \tA = \tA X$. Then  $M= \tA A^{-1} L A \tA^{-1}$. Let $U=A \tA^{-1}$. Note that  $LU=UM$.  Using the identity (3.4)  it is easy to obtain 
$U=L D_q -q D_q  M.$ Combining these equations we get 
$$    L^2 D_q - (q+1) L D_q M + q D_q M^2 =0.    \eqno(3.6)$$

We define the {\it extended Hahn's class} of $\cD_q$, denoted by $\cJ_{q,+}$, as the set of all  MOPs $A$ such that its matrix $L \in \cT$ of recurrence coefficients satisfies the quadratic matrix equation 
$$    L^2 D_q - (q+1) L D_q M + q D_q M^2 = t D_q,       \eqno(3.7)$$
for some $M \in \cT$ and some complex number $t$. 

Let
$$L=\left[  \begin{matrix} \beta_0 & 1 & 0 & 0 & \ldots \cr \alpha_1 & \beta_1 & 1 & 0 & \ldots \cr 0 & \alpha_2 & \beta_2 & 1 & \ldots \cr \vdots &\vdots & \vdots & \vdots & \ddots \cr \end{matrix} \right], \qquad 
	M=\left[  \begin{matrix} \tilde \beta_0 & 1 & 0 & 0 & \ldots \cr \tilde \alpha_1 & \tilde \beta_1 & 1 & 0 & \ldots \cr 0 & \tilde \alpha_2 & \tilde \beta_2 & 1 & \ldots \cr \vdots &\vdots & \vdots & \vdots & \ddots \cr \end{matrix} \right].
	\eqno(3.8)$$
Let  $\sigma_k= \beta_0+\beta_1+ \cdots + \beta_k$ and $\tilde \sigma_k =  \tilde \beta_0+ \tilde \beta_1+ \cdots + \tilde \beta_k$ for $k \geqslant 0$.  We write  the coefficients $\beta_k$ and $\tilde \beta_k$  in terms of the $\sigma_k$ and $\tilde \sigma_k$  and then we solve the scalar equations obtained by comparing corresponding entries in both sides of (3.7). 
We obtain first 
$$   \tilde \sigma_k = \frac{[k+1]}{ [k+2]} \sigma_{k+1}, \qquad  k \geqslant 0,    \eqno(3.9)$$ 
and 
$$   \tilde \alpha_k = \frac{[k] g(k+1)}{q \, [k+1] g(k) } \alpha_{k+1},
\qquad k \geqslant 1, \eqno(3.10)$$
where the function $g$ is defined by
$$ g(k) = [k-2] c_0 +  c_1,  \qquad k \geqslant 0, \eqno(3.11)$$
and   the constants   $c_0$ and  $c_1$ are defined by  
$$ c_0 = [3] (q \sigma_0 -\sigma_1) + \sigma_2, \qquad c_1 = [3] (\sigma_0 + \sigma_1 -\sigma_2).   \eqno(3.12)$$

From the rest of the scalar  equations obtained from (3.7), we obtain 
$$  \sigma_k=   \frac{[k+1] ( [k] q^{-1} c_2+ \sigma_0 g(1)) }{g(2k+1)}  ,\qquad k \geqslant 0, \eqno(3.13)$$
and, for $k \geqslant 1$ 
$$	\alpha_k=  
\frac{q\, [k]\, g(k-1)}{g(2k)g(2k-2)} 
\left( q^{k-2} g(2) \alpha_1  
  + [k-1] g(k)  W(k) \right), \eqno(3.14)$$
where 
$$W(k)= 
  \frac{ q^{k-3} ((q^{k-1} +1)  \sigma_0 c_0 -c_2)
        ((\sigma_0 c_0-c_2) q^{k-2} +\sigma_0 (c_0- (q-1) c_1))
    }{g(2k -1)^2  }-t,  $$ 
and
the auxiliary constant $c_2$ is defined by
$$ c_2= [4] \sigma_0 \sigma_2 -[3] \sigma_0 \sigma_1 -q^2 \sigma_1 \sigma_2. \eqno(3.15)$$ 

Note that the sequence  $\sigma_k$ is independent of  $t$ and that $\alpha_k$ is a rational function of $q^k$ with numerator and denominator of degree less than or equal to  8. Both sequences, $\sigma_k$ and $\alpha_k$, are completely determined by the parameters $\alpha_1, \beta_0, \beta_1, \beta_2,$ and $t$.  

The expression for $\alpha_k$ is obtained by using Maple to compute $\alpha_k$ for $k=2,3,\ldots, 10$ and then  using polynomial interpolation to find the numerator in (3.14) as a function of $q^k$. Formula (3.14) is a minor modification of the corresponding formula in \cite{Rec} due to a change in the definition of the  auxiliary constants.
It is easy to obtain the limits of $\sigma_k$ and $\alpha_k$ as $q$ goes to 1. The resulting expressions coincide with the formulas for the recurrence coefficients of the classical orthogonal polynomial sequences obtained in \cite{Mops}. 

Some additional properties of the elements of the extended $q$-Hahn class can be found in \cite{Rec}.

\section{Zeroes and poles of the recurrence coefficients}

In this section we obtain formulas for the coefficients $\sigma_k$ and $\alpha_k$ in which the  numerators and denominators are  completely factorized as products of polynomials of degree one in $q^k$. We also obtain some properties of the zeroes and poles of  the sequence $\alpha_k$.

\begin{theorem} Let      $y,d_1,d_2,d_3$  be complex numbers. Define the rational functions

$$Z(x)=\frac{(x-1) (x-y)}{(x^2-q y) (x^2 -y)^2 (q x^2 -y)} \prod_{j=1}^3( (x-d_j) (d_j x -y)),\eqno(4.1) $$
and 
$$V(x)=\frac{(q x -1) (k_1 q x -k_2)}{(q-1) (q^2 x^2-y)}, \eqno(4.2)$$
 where  the constants $k_1$ an $k_2$ are defined by  $$k_1=y+ d_1 d_2 + d_1 d_3 + d_2 d_3, \qquad  k_2=y (d_1+d_2+d_3) +d_1 d_2 d_3.$$
Suppose that  $V(1)\ne 0$ and  $Z(q^k)\ne 0 $ for $k \geqslant 1$. Define the sequences 
$$\alpha_k= \frac{ \alpha_1  Z(q^k)}{Z(q)}, \qquad  k \geqslant 1, \eqno(4.3)   $$
and   
$$\sigma_k=\frac{\sigma_0 V( q^k )}{V(1)} , \qquad k \geqslant 0,  \eqno(4.4) $$
where $\alpha_1$ is a nonzero number and $\sigma_0$ is defined by
$$\sigma_0^2= \frac{ \alpha_1  V(1)^2}{     Z(q)}. \eqno(4.5)$$
Let $L$ be the matrix determined by the sequences $\alpha_k$  and $\sigma_k$ as in (3.8), with $\beta_k=\sigma_k-\sigma_{k-1}$ and let $A$ be the unique matrix in $\cG$ that satisfies $LA=AX$. Then $A$ is an element of  the extended $q$-Hahn class $\cJ_{q,+}$ that satisfies (3.7) with
$$t=- \frac{\alpha_1   (q-1)^2 d_1 d_2 d_3}{q^2 Z(q)}, \eqno(4.6) $$ and $M$ defined by equations (3.8), (3.9) and (3.10).
\end{theorem}

The proof of this theorem is a straightforward computation that verifies that the quadratic matrix equation (3.7) is satisfied. This is easily done using Maple. 
We will explain how we found the functions $Z(x)$ and $V(x)$ later. 
Observe that the parameters $d_1,d_2,d_3$ play a symmetric role in the construction of the recurrence coefficients. 

In order to simplify the notation from now on  we write $ y=p^2$   and $q=r^2$, that is $\pm p$  are the square roots of $y$ and $\pm r$ are the square roots of $q$.

From (4.1) we see that the roots of $Z(x)$ are  $ 1, p^2, d_1,d_2,d_3,p^2/d_1,p^2/d_2,p^2/d_3,$ and its poles are $ p,-p,p,-p,pr,-pr, p/r, -p/r. $
Observe that if $p\ne 0$ then  the set of zeroes and the multi-set of poles of $Z(x)$ are invariant under the map that sends $z$ to $p^2/z$  for $z$ in the extended complex plane.
Because of this fact, if $p, d_1,d_2,d_3$ satisfy some conditions then the sequence $\sigma_k$ of the previous theorem can be replaced by a different sequence $\hat\sigma_k$ and the matrix associated with $\alpha_k$ and $\hat\sigma_k$ is also a MOP that belongs to the extended $q$-Hahn class. This is expressed in detail in the following theorem.

\begin{theorem} With the definitions and notation of the previous theorem, 
suppose now  that $p, d_1,d_2,d_3$ satisfy  $p\ne 0$ and   $k_2 q - y k_1\ne 0$.  
Define 
$$V_a(x)= \frac{(q x -1) (k_2 q x - y k_1)}{(q-1) (q^2 x^2-y)}, \eqno(4.7)$$
and
$$\hat\sigma_k= \frac{\hat\sigma_0 V_a( q^k )}{V_a(1)} , \qquad k \geqslant 0,  \eqno(4.8) $$
where 
$$\hat\sigma_0^2=\frac{\alpha_1 V_a(1)^2}{y Z(q)}. \eqno(4.9)$$
Then the MOP $\hat A$ determined by the recurrence coefficients $\alpha_k$ and $\hat\sigma_k$ 
is in the  extended $q$-Hahn class $\cJ_{q,+}$ and  satisfies (3.7) with  $t$ defined by (4.6).
Furthermore, $\hat\sigma_k$ coincides with $\sigma_k$ if and only if $p=\pm d_j$ for some $j$ such that $1 \leqslant j \leqslant 3$. 
\end{theorem} 

The proof is a direct verification of equation (3.7), as in the previous theorem. 

From Theorem 4.2 we see  that, in many cases,  for each $(y,d_1,d_2,d_3)$ there are two non-equivalent 
 $q$-orthogonal polynomial sequences that correspond to two distinct matrices in $\cJ_{q,+}$.

The degrees of the numerator and the denominator of the function $Z(x)$ depend on how the multisets of zeroes and poles are related to each other and also on whether some of the parameters $(y,d_1,d_2,d_3)$ are equal to zero. We will write  $\deg(Z)=(n,m)$ if $n$ is the degree of the numerator and $m$ is the degree of the denominator. It is clear that  $n$ and $m$ can take values between 0 and 8, but not all pairs are possible because cancellation of zeroes and poles  occurs by pairs.    
Let us note that $\alpha_k$ is a rational function of $q^k$ which has the same degrees of the numerator and denominator as $Z(x)$.

If $f$ is a function of $q$ we define the $q$-reversed function $R_q(f)$ by 
 $R_q(f)(q)=f(1/q)$.  
The map $R_q$ can be used to transform matrices. The following proposition clearly holds.

\begin{prop}  Suppose that $A$ is the  element of $\cJ_{q,+}$ that corresponds to the sequences 
 $\alpha_k$ and $ \sigma_k$ and the parameters $ \sigma_0$ and $t$, defined as in Theorem 4.1. Then the matrix $R_q(A)$ that corresponds to the sequences $R_q(\alpha_k)$ and $R_q(\sigma_k)$ and the parameters $R_q( \sigma_0)$ and $R_q(t)$ is also an element of $\cJ_{q,+}$. 
\end{prop}

Note that  $\deg(R_q(\alpha_k))$ is in general not equal to $\deg(\alpha_k)$.

It is easy to verify that if $y=1$ then $\alpha_k$ and $t$ are invariant under $R_q$ and  $R_q(\sigma_k) = \hat\sigma_k$, and  $ R_q( \sigma_0^2) = \hat \sigma_0^2.$ 

Now we present a brief description of the procedure used to obtain the factored formula (4.3) for the sequence  $\alpha_k$ from (3.14). We replace first  $q^k$ with $x$ in (3.14). The factors $[k]$ and $g(k-1)$ yield polynomials of degree one in $x$. The first one gives $x-1$ and the second one $x-y$, where $y$ is expressed in terms of the constants $c_0$ and $c_1$ defined in (3.12). The denominator of (3.14) is already factored and its roots are expressed in terms of $y$. The remaining factor in the numerator is a polynomial of degree six. 
Suppose that $c_0 \ne 0$ and write $t= w \alpha_1$, $c_3=u c_1$ and $\sigma_0^2=v \alpha_1$. These substitutions give us a polynomial in $x$ of degree six that we call $F(x)$.   Now solve $F(d_1)=0$ for the parameter $w$ and substitute the result in $F(x)$.  This gives us a new polynomial that factors as $(x-d_1)(d_1 x -y) F_1(x)$ where $F_1$ is a polynomial of degree four in $x$ that depends also on $v$ and $u$. Then we solve $F_1(d_2)=0$ for $v$ and substitute the result in $F_1$. We get a polynomial that factors as $(x-d_2)(d_2 x -y)F_2(x)$, where $F_2(x) $ is a quadratic polynomial that factors as $F_2(x)=G_1(x)G_2(x)$. Solving $G_1(d_3)=0$ for $u$ we obtain a solution $u_1$ that depends on $y, d_1,d_2,d_3$. Substitution of $u=u_1$ in $F_2$ gives $F_2(x)=K (x-d_3)(d_3 x -y)$ and this completes the factorization of $F(x)$. If instead of solving $G_1(d_3)=0$  for $u$ we solve $G_2(d_3)=0$ for $u$ we obtain a different solution $u_2$.  This is what produces two possible different formulas for $\sigma_k$. Observe that what makes possible the symbolic factorization of $F(x)$ is that its roots appear by pairs, $d_j$ and $y/d_j$.  

For the case $c_0=0$ we found that the formulas of Theorem 4.1 also hold. For that case the parameters $y, d_1,d_2,d_3$ must satisfy certain relation.

\section{The extended $q$-Hahn class and the $q$-Askey scheme}

In this section we show how to obtain some elements of the $q$-Askey scheme using our general representation of the elements of the extended $q$-Hahn class in terms of the parameters $y, d_1,d_2,d_3$  and the arbitrary nonzero scaling parameter $\alpha_1$. 
We have verified that all the examples in \cite[Chapter 14]{Hyp} can be obtained in such way.
Recall that $\sigma_k=\beta_0+\beta_1+\cdots+\beta_k$  for $k \geqslant 0$.

The recurrence coefficients for the monic Askey-Wilson polynomials  (\cite[14.1.5]{Hyp})  are obtained by taking   

$$y=\frac{a b c d}{q^2}, \ \ \ \   d_1=\frac{bc}{q},\ \  \ \ d_2=\frac{cd}{q},\ \  \ \ d_3=\frac{bd}{q}. \eqno(5.1)$$
These equations can be solved for the Askey-Wilson parameters. We obtain
$$a=y \sqrt{\frac{q}{d_1 d_2 d_3}}, \ \ \  b=\sqrt{\frac{q d_1 d_3}{d_2}},\  \ \  c=\sqrt{\frac{q d_1 d_2}{d_3}},\  \ \ d=\sqrt{\frac{q d_2 d_3}{d_1}}. \eqno(5.2)$$
If some $d_j=0$ then the Askey-Wilson parameters are not defined.

For the  monic $q$-Racah polynomials the recurrence coefficients (\cite[14.2.4]{Hyp}) are obtained with
$$ y= \frac{1}{\alpha \beta}, \ \ \ d_1=\frac{1}{\alpha}, \ \ \  d_2=\frac{1}{\gamma},\ \ \ d_3=\frac{\delta}{\alpha}. \eqno(5.3)$$
For the Askey-Wilson and the $q$-Racah polynomials we have $\deg(\alpha_k)=(8,8)$.

Recall that $q=r^2$ and $y=p^2$. 
Taking $y=1, d_1=-1,d_2=r, d_3=-r$ we obtain the constant sequences of  coefficients $\alpha_k=\alpha_1$, 
and   $\sigma_k=\sigma_0$  with $\sigma_0^2=\alpha_1$ and $t= -\alpha_1 (q-1)^2 /q$.   
Taking  $y=q,d_1=r, d_2=-r, d_3=-1$ we obtain $\alpha_k=\alpha_1 /2$ for $k\geqslant 2$ and $\sigma_k=0$ for $k \geqslant 0$, with $t=-\alpha_1 (q-1)^2/2 q$. These two examples show that there are  sequences with  $\deg(\alpha_k)=(0,0)$ in the extended $q$-Hahn class.

In the case  $y=1$ and $d_1,d_2,d_3$  arbitrary we obtain  $\deg(\alpha_k)=(6,6)$ and the reversed sequence also has degrees (6,6). This family of recurrence coefficients is not mentioned in \cite{Hyp} and therefore it may be new. In fact, none of the 32  examples in  \cite[Chapter 14]{Hyp} has  $\deg(\alpha_k)=(6,6)$.

If we take $y=1, d_1=-1$ and $d_2,d_3$ arbitrary, we obtain a family of sequences with  $\deg(\alpha_k)=(4,4)$ and such that the reversed sequence has also degrees (4,4). This family is also possibly new.  The  element of this family that has $d_2=d_3=1$  gives us 
$$\alpha_k=  \frac{[3]  \alpha_1 [k]^4}{[2k+1] [2k-1]},\  \ \ \sigma_k=0,\  \ \ t=-\frac{[3]\alpha_1}{q}. \eqno(5.4)$$  

Another interesting example is obtained  when  $y=1, d_1=d_2=d_3=0. $ For this case we have
$$\alpha_k=\frac{[3] (q+1)^2 q^{3(k-1)} \alpha_1}{[2k+1][2k-1] (q^k+1)^2 }, \ \ \sigma_k=\frac{q^k (q+1) \sigma_0}{q^{k+1}+1}, \ \ \sigma_0^2=-  [3] \alpha_1 , \ \ t=0, \eqno(5.5)$$
and  $R_q(\alpha_k)=\alpha_k$, $R_q(\sigma_k)=q^k \sigma_k$, and  $R_q(\sigma_0^2)=\sigma_0^2/q^2$.
 Since $t=0$ the corresponding MOPS  $A$ and $R_q(A)$ are  in the $q$-Hahn class. Observe that  $\deg(\alpha_k)=(3,6)$ and  if $\alpha_1 >0$ and  $q>0$  then $\sigma_0^2$ must be negative. 

 The map $q \rightarrow q^{-1} $ induces a pairing of matrices of recurrence coeficients. In order to deal with this pairing,  and with the zeroes and poles at infinity of the rational functions $Z$ and $V$,  it is convenient to introduce the following definitions.
$$Z^*(x)= R_q(Z(1/x))=    
\frac{q^2 (x-1) (y x-1 )}  {(y x^2-q ) (y x^2 -1)^2 (y q x^2 -1)} \prod_{j=1}^3( (y x-d_j) (d_j x -1)),\eqno(5.6) $$
$$V^*(x)=R_q(V(1/x))= \frac{q (q x -1) (k_2 q x -k_1)}{(q-1)(y q^2 x^2-1)}, \eqno(5.7)$$
$$V_a^*(x)=R_q(V_a(1/x))=  \frac{q (q x -1) (y k_1 q x -k_2)}{(q-1)(y q^2 x^2-1)}, \eqno(5.8)$$
$$(\sigma_0^*)^2=\frac{\alpha_1 V^*(1)^2}{Z^*(q)}, \eqno(5.9)$$
and 
$$(\hat\sigma_0^*)^2=\frac{\alpha_1 V_a^*(1)^2}{y   Z^*(q)}. \eqno(5.10)$$

It is easy to verify that 
$$\alpha_k^* =R_q(\alpha_k)=\frac{\alpha_1 Z^*(q^k)}{Z^*(q)}, \qquad k \geqslant 1, \eqno(5.11)$$
$$\sigma_k^* =R_q(\sigma_k)= \frac{\sigma_0^* V^*(q^k)}{V^*(1)}, \qquad k \geqslant 0, \eqno(5.12)$$   
$$\hat\sigma_k^* = R_q(\hat\sigma_k)=  \frac{\hat\sigma_0^* V_a^*(q^k)}{V_a^*(1)}, \qquad k \geqslant 0, \eqno(5.13) $$ and
$$t^*= R_q(t)=\frac{-\alpha_1 d_1 d_2 d_3 (q-1)^2}{Z^*(q)}.
 \eqno(5.14)$$
These last four formulas allow us to compute in a simple way the recurrence coefficients when some of the parameters $p, d_1,d_2,d_3$, or their reciprocals, are equal to zero. There is no need to deal with limits.
 For example, if we take $p=0, d[2]=d[3]=0$ and $d[1]=1/v$ in (5.11)  and (5.12) and then we put $v=0$ we obtain the recurrence coefficients of the discrete $q$-Hermite polynomials.

There are cases in which the relation between $\sigma_0$ and $\alpha_1$ determined in equations (4.5) and (4.9) is not a necessary condition and thus $\sigma_0$ becomes an arbitrary parameter.

\section{Final remarks}
In this paper we have only considered properties of the elements of the extended $q$-Hahn class that are directly related with the recurrence coefficients.
There are many other properties, such as the generating functions, weight functions, Rodrigues formulas, and $q$-difference equations  that we have not studied yet. 

In \cite{Rec} we proved that the elements of the extended $q$-Hahn class, expressed in terms of certain polynomial basis, satisfy a second order $q$-difference equation of Bochner type. We also found explicit formulas for the coefficients of the $q$-difference operator.  
One important topic is the study of the connections between the theory of $q$-difference operators on lattices and the tridiagonal matrices of change of bases that map  the extended $q$-Hahn class into the classical $q$-Hahn class, or equivalently, that transform the classical $q$-difference operators into operators that, in some cases, may be interpreted as difference operators on lattices.

For the extended  Hahn class of discrete orthogonal polynomials we have obtained results, analogous to the ones reported here, that will appear elsewhere.

\section{Aknowledgements}
This research was partially supported by grant 220603 from CONACYT, Mexico.


\begin{thebibliography}{00} 
\bibitem{Alv} R. \'Alvarez-Nodarse, R.  Sevinik Adig\"uzel, and H. Taseli, On the orthogonality of $q$-classical polynomials of the Hahn class, SIGMA Symmetry Integrability Geom. Methods Appl. 8 (2012), Paper 042, 30 pp.
\bibitem{H} W. Hahn, \"Uber die Jacobischen Polynome und zwei verwandte Polynomklassen, Math. Zeit., 39 (1935) 634--638.

\bibitem{Hyp} R.  Koekoek, P. A. Lesky, and R. F. Swarttouw, Hypergeometric orthogonal polynomials and their q-analogues,  Springer-Verlag,  Berlin, Heidelberg, 2010.
\bibitem{KS} W. Koepf  and D.  Schmersau, 	Recurrence equations and their classical orthogonal polynomial solutions, Orthogonal systems and applications, 
	Appl. Math. Comput. 128 (2002)  303--327.
\bibitem{Koor} T. H. Koornwinder,  The Askey scheme as a four-manifold with corners. Ramanujan J., 20 (2009)  409--439.

\bibitem{MBP} F. Marcell\'an, A. Branquinho, and J. Petronilho, Classical orthogonal polynomials: a functional approach, Acta Appl. Math.,  34 (1994) 283--303.	
\bibitem{Med} J. C. Medem, R. \'Alvarez-Nodarse, F. Marcell\'an, On the $q$-polynomials: a distributional study, J. Comp. Appl. Math. 135 (2001) 157--196.
\bibitem{Mops} L. Verde-Star,	Characterization and construction of classical orthogonal polynomials using a matrix approach, Linear Algebra  Appl. 438 (2013) 3635--3648.
\bibitem{Rec} L. Verde-Star, Recurrence coefficients and difference equations of classical discrete and $q$-orthogonal polynomial sequences,     Linear Algebra  Appl. 440 (2014) 293--306. 

\end{thebibliography}
\end{document}